\documentclass[11pt,a4paper]{article}
\usepackage{amsthm,amsmath,amssymb}
\usepackage{amsfonts}
\usepackage{tikz}
\usepackage{relsize}
\usepackage{verbatim}

\newtheorem{thm}{Theorem}[subsection]
\newtheorem{prop}[thm]{Proposition}

\theoremstyle{definition}

\newtheorem{lem}[thm]{Lemma}

\theoremstyle{remark}

\def\ket[#1]{\left| #1 \right>}
\def\bra[#1]{\left< #1 \right|}

\def\cat{\mathcal C}
\def\A{Abelian }
\newcommand{\beq}{\begin{equation}}
\newcommand{\eeq}{\end{equation}}
\newcommand{\bc}{\begin{center}}
\newcommand{\ec}{\end{center}}

\newif\ifPrivateMode
\PrivateModefalse
\long\def\/*#1*/{}

\title{From tensor category to Temperley-Lieb algebra representation}
\author{Peter E. Finch, Zolt\'an K\'ad\'ar, Paul Martin
\\\\School of Mathematics, University of
  Leeds, Leeds, LS2 9JT, United Kingdom
}
\date{\today}

\begin{document}
\newcommand{\epsrc}{This work was supported by EPSRC under grant EP/I038683/1}
\maketitle


\begin{abstract}
We construct a representation of the Temperley-Lieb algebra from a multiplicity-free semisimple monoidal Abelian category, with two simple objects 
$\lambda$ and $\nu$ such that $\lambda\otimes\nu$ is simple and $Hom_{\cal C}(\lambda\otimes\lambda, \nu)$ is not empty. 
A self-contained manual to tensor categories is also provided as well as a summary of the best known example of the construction:
Schur-Weyl duality for $U_q(sl_2))$.    
\end{abstract}
\section{Introduction}
The Temperley-Lieb algebra $TL_{n}(\delta)$ has a faithful representation $\theta$ on tensor space as linear maps in 
End$\left(V^{\otimes n}\right)$, where
$V$ is the fundamental two dimensional module of $U_q(sl(2))$ for almost all values of $\delta\in {\mathbb C}$ \cite{CFS,PH}. 
In this representation the generator $U_i$ is mapped to the map, which is the identity on all tensor factors except for $i,i+1$,
where it is the composition in End$(V\otimes V)$ of the projection to the 1-dimensional module $V^{(0)}$ of $U_q(sl(2))$ with the
injection to $V\otimes V$ \cite{CFS,ZTQC}. This enables to identify Temperley-Lieb algebra elements with Homs in $U_q(sl(2))$. 
In other words, we have the semisimple monoidal, Abelian category Rep\,$U_q(sl(2))$ containing simple objects $V$ and $V^{(0)}$, 
and we can build a representation of $TL_n(\delta)$. 

%
%

One can realise that the proof is independent of the details of the category that we started from and gives rise to a general theorem. It states that in defining
$\theta$ both the fundamental module $V$ and $V^{(0)}$ ``can be substituted'' with some simple objects $\lambda$ and $\nu$ 
along with the map $\cup\circ\cap$ in the non-trivial part of $\theta(U_i)$ with projections in End$(\lambda\otimes \lambda)$ 
without destroying the homomorphism property.
	
The context, where this result first came from was integrable models. If $H$ is a quasi-triangular Hopf algebra, one may wish to determine 
the centraliser algebra of \{End$(\lambda(a)^{\otimes n}), a\in H\}$, and determine whether it has the Temperley-Lieb, Hecke or 
Birman-Murakami-Wenzl algebra as a quotient. In case the answer is yes, a solution of the Yang-Baxter equation can be found. This happens
in Interacting Round a Face models \cite{Pasquier1988} and it pops up in the correspondence between quantum spin and anyon chains as well
\cite{Finch1,Finch2,FTLTKWF,GATHLTW}.

Before stating the theorem we give an introduction to tensor categories, monoidal systems and the diagram 
calculus encoding the category composition in several fields of physics and mathematics: e.g., comformal field theory, 
topological quantum computation. Then the theorem will be stated and proved in section \ref{prdr} using 
the machinery of diagram calculus. For an example, in section \ref{SW}
 we will explain the known case of $U_q(sl(2))$ and show how the pairs of sequences basis \cite{Paulbook} of 
$TL_n(\delta)$ is presented in terms of diagrams. Further questions and examples are mentioned in the last section. In the appendix, the 
matrix elements of the projections are given in terms of F-symbols, which establishes equivalence with the proof given in \cite{Finch3} in
terms of monoidal systems.
\section{Definitions}
Here we provide the definitions of the notions needed to define semisimple monoidal abelian tensor categories. 
Following the terminology from \cite{yama} we will call these tensor categories. The standard references are \cite{BK,MacLane}. 
Then we introduce the diagram calculus and establish the equivalence between tensor categories and the so-called monoidal
systems, a set of data containing all information of the category.  
\subsection{Some basic notions from categories}
This small section contains some general definitions we need later. 
\begin{itemize}
\item An arrow $i \in $Hom$_\cat(a,b)$ is an {\em injection} (a.k.a. coproduct injection) if for any arrow 
$f\in$ Hom$_\cat(a,c)\;\exists!\;h\in$ Hom$_\cat(b,c)$ such that $f=h\circ i$. 
\item An object $t$ in a category $\cat$ is {\em terminal} 
if $|$Hom$_\cat(a,t)|=1\;\forall a\in$Obj$_\cat$, that is, there is precisely one arrow in the category from every object to $t$.
\item An object $s$ in a category $\cat$ is {\em initial} 
if $|$Hom$_\cat(s,a)|=1\;\forall a\in$Obj$_\cat$.
\item An object $z$ in a category $\cat$ is {\em zero} if it is both initial and terminal.
\item In a category $\cat$ with zero object, the {\em kernel} of an arrow $f\in$ Hom$_\cat(a,b)$ is an arrow $k\in$ Hom$_\cat(c,a)$ such that 
$f\circ k=0$ and $\forall h\in$ Hom$_\cat(c',a)$ with $f\circ h=0, \exists!\;l\in $Hom$_\cat(c',c)$ suct that $k\circ l=h$. Note that 
in general, $0\in $Hom$_\cat(d,e)$ is the unique composition $d\to z\to e$, where $z$ is the zero object in $\cat$. 
\item In a category $\cat$ with zero object, the {\em cokernel} of an arrow $f\in$ Hom$_\cat(a,b)$ is an arrow $u\in$ Hom$(b,c)$ such that
$u\circ f=0$ and $\forall h\in$ Hom$_C(b,c')$ with $h\circ f=0, \exists!\;l\in $Hom$_\cat(c,c')$ suct that $l\circ u=h$.
\item The object a in a category $\cat$ is a {\em product} of the objects $a_1,a_2$ if $\exists \pi_i:a\to a_i, i=1,2$ such that for any 
object $b$ and arrow $f_i:b\to a_i, \exists! f:b\to a$ with $f_i=\pi_i\circ f, i=1,2$.  
\item The object a in a category $\cat$ is a {\em coproduct} of the objects $a_1,a_2$ if $\exists \iota_i:a_i\to a, i=1,2$ such that for any
object $b$ and arrow $f_i:a_i\to b, \exists! f:a\to b$ with $f_i=f\circ \iota_i, i=1,2$.
\item A {\em functorial morphism} for two functors $F,G:\cat\to\cat'$ is a collection of arrows $\phi_a:F(a)\to G(a), 
a\in$Obj$(\cat)$ such that for every $f\in $Hom$_\cat(U,V): \phi_V\circ F(f)=G(f)\circ \phi_U$.
\end{itemize}
\subsubsection{Additive categories}
We would like to be able to add morphisms. In an {\em Ab} category each Hom-set is an additive abelian group and the composition is bilinear:
\[(g+g')\circ (f+f')=g\circ f+g\circ f'+g'\circ f+g'\circ f'\] 
for any arrows $g,g'\in$ Hom$_\cat(b,c), f,f'\in $Hom$_\cat(a,b)$ and $+$ is the addition in the corresponding abelian groups. 
The object $a$ in an Ab-category $\cat$ is a {\em biproduct} of the objects $a_1,a_2$, if there are arrows shown in the 
diagram:\vspace{0.2cm}

\[a_1\substack{\xleftarrow{\pi_1} \\[-0.7ex] \xrightarrow[\iota_1]{}} a\substack{\xrightarrow{\pi_2} 
\\[-0.7ex] \xleftarrow[\iota_2]{}}a_2\;\;\mbox{with}\;\;\pi_i\circ \iota_i=\mbox{id}_{a_i}, i=1,2\;\;\mbox{and}\;\; 
i_1\circ \pi_1+i_2\circ \pi_2=\mbox{id}_a.\]

\hspace{-0.6cm}Note that if there are binary biproducts, then iteration for given $a_1, a_2,\dots, a_n$ yields a 
biproduct $a\equiv\oplus_i^n a_i$ characterized (up to isomorphism in $\cat$) by
\[a_j\xrightarrow{\iota_j}\bigoplus_i a_i\xrightarrow{\pi_k}a_k,\;\;j,k=1..n,\;\;
\mbox{with}\;\;\sum_i^n\iota_i\circ\pi_i=\mbox{id}_a,\;\;\pi_j\circ\iota_k=\mbox{id}_{a_j}\delta_{j,k}\]
Note that each Hom set has a zero arrow (the unit of the addition in the abelian group) 
even if the category has no zero object in the above sense.  
\vspace{0.2cm}

\hspace{-0.6cm}The category $\cat$ is an {\em additive} category if it is an Ab category with
\begin{itemize}
\item zero object
\item binary biproducts between any two objects
\end{itemize}
Note that there are biproducts in an Ab-category iff there are products iff there are coproducts \cite{MacLane}. 
We will call the biproducts direct sum (as is the usual notation in module categories or the category of 
abelian groups). 
\subsubsection{Abelian categories}
The additive category $\cat$ is an {\em Abelian} category if 
\begin{itemize}
\item there exists kernel and cokernel for every arrow and
\item every monomorphisms is a kernel and every epimorphism is a cokernel
\end{itemize}
Note, that generally every kernel is a monomorphism and every cokernel is an epimorphism, but the converse is not always 
true. The proofs of these statements are e.g. in \cite{MacLane}.
\prop{\cite{MacLane} In Abelian category every arrow has a factorisation $f=ker(coker f)\circ coker(ker f)$}. 
\subsection{Monoidal tensor categories}
A category $\cat$ is {\em monoidal} if it has the following data (additional to just being a category):
\begin{itemize}
\item a bifunctor $\otimes: \cat \times \cat  \to \cat$
\item a functorial isomorphism (associativity) $\alpha_{a,b,c}:(a\otimes b)\otimes c\to a\otimes (b\otimes c)$ of
functors $\cat\times\cat\times\cat\to\cat$
\item a unit object ${\bf 1}$ and functorial isomorphisms $\lambda_a:{\bf 1}\otimes a\to a, \rho_a: a\otimes {\bf 1}\to a$ 
for $a\in$Obj$(\cat)$.
They have to satisfy the asociativity axiom:
\item If two expressions $X_1$ and $X_2$ are obtained from $a_1\otimes a_2\otimes \dots \otimes a_n$ by inserting brackets 
and ${\bf 1}$'s, then all isomorphism composed of $\alpha,\lambda,\rho$ and their inverses have to be equal. 
\end{itemize}
\thm{(MacLane Coherence Theorem). Suppose we are given the data $(\cat,\otimes,\alpha,\lambda,\rho)$. $\cat$ is a monoidal category
iff the following classes of relations are satisfied.
\begin{itemize}
\item The pentagon relations are the equalities of the following two compositions for any quadruple $(a,b,c,d)$ of objects in $\cat$:
\[\begin{array}{c}
((a\otimes b)\otimes c)\otimes d\xrightarrow{\alpha_{a,\,b,\,c}\otimes\mbox{{\small id}}_d}(a\otimes (b\otimes c))\otimes d
\xrightarrow{\alpha_{a,\,b\otimes c,\,d}}a\otimes((b\otimes c)\otimes d)\xrightarrow{\mbox{{\small id}}_a\otimes\alpha_{b,\,c,\,d}}
a\otimes (b\otimes (c\otimes d))\\=\\
((a\otimes b)\otimes c)\otimes d\xrightarrow{\alpha_{a\otimes b,\,c,\,d}}(a\otimes b)\otimes (c\otimes d)
\xrightarrow{\alpha_{a,\,b,\,c\otimes d}}a\otimes (b\otimes (c\otimes d))
\end{array}\] 
\item
The triangle relations are the equalities of the following compositions for any pair $(a,b)$ of objects in $\cat$:
\[\begin{array}{c}(a\otimes{\bf 1})\otimes b\xrightarrow{\alpha_{a,\,{\bf 1},\,b}}a\otimes ({\bf 1}\otimes b)
\xrightarrow{{\small\mbox{id}_a}\otimes \lambda_b}a\otimes b\\=\\
(a\otimes{\bf 1})\otimes b\xrightarrow{\rho_a\otimes {\small\mbox{id}_b}}a\otimes b
\end{array}\]  
\end{itemize}
 } 
\proof{In \cite{MacLane,BK}.} 
\subsection{Semisimple categories \label{sc}}
An object $a$ in an \A category $\cat$ is called {\em simple} \cite{BK}
\begin{itemize} 
\item if it is not isomorphic to the zero object and
\item if any injection in Hom$_\cat(b,a)$ is either $0$ or isomorphism for any object $b$. 
\end{itemize}
Let $I$ denote the set of equivalence classes of simple objects and 
choose a representative simple object $\lambda_i$ from each class $i\in I$.  
An \A category is {\em semisimple} if it has countably many isomorphism classes of simple objects and 
any object $a$ is isomorphic to a direct sum of simple ones $a\cong \oplus_{\lambda\in \hat{I}_a}\lambda$ and 
the multiset $\hat{I}_a$ of simple objects is finite. Let $N_a^\lambda$ denote the number of appearance of $\lambda$ in $\hat{I}_a$.  
Example: Rep$\,U(sl(2))$ the representation category of the algebra $U(sl(2))$. Simples are indexed by integers and the fusion 
rules are given by 
\[i\otimes j\cong \bigoplus_{k=0}^{(i+j-|i-j|)/2} |i-j|+2k\]
Some concrete instances: 
\[\begin{array}{c} 1\otimes 1\cong 0\oplus 2 \\\\
a\equiv (1\otimes 1)\otimes 1\cong (0\oplus 2)\otimes 1\cong (0\otimes 1)\oplus (2\otimes 1)\cong 1\oplus 1\oplus 3
\end{array}\ .\]
Using the notation above $N^{1}_a=2,\,N^{3}_a=1,\,N_a^\lambda=0$ for 
any integer lambda other than $1$ or $3$ and $\hat{I}_a=[1,1,3]$. Note that the monoidal unit here is denote by $0$.

Let $k$ be a field of characteristic zero. Following the terminology from \cite{yama} 
we define {\em tensor} category $\cat$ to be a semisimple \A monoidal category, where the tensor product distributes
over the addition of morphisms with the following additional properties:
\begin{itemize}
\item All Hom spaces are k-vector spaces,
\item the arrow composition is $k$-bilinear and
\item End$(\lambda_i)=k$ for all $i\in I$ and Hom$(\lambda_i,\lambda_j)=0, i\neq j$. 
\end{itemize}
Note that the last property is automatically fulfilled in case $k$ is algebraically closed (Schur's lemma). From now on we always
work with such a category unless stated otherwise.
Let us choose a representative of all ismomorphism classes of simple objects and let the letter starting from $\lambda$ in 
the Greek alphabet denote these representatives from now on. Let us introduce the notation 
$N_{\mu,\nu}^\lambda\equiv N_{\mu\otimes \nu}^\lambda$ and
 $\hat{I}_{\mu,\nu}\equiv\hat{I}_{\mu\otimes \nu}$.
Due to the axioms (semisimplicity and existence of all biproducts) for all pairs $(\mu,\nu)$ of labels of simple objects 
there exist a collection of arrows
\begin{equation}\mu\otimes \nu\xrightarrow{\pi^{\lambda,\alpha}_{\mu,\nu}}\lambda
\xrightarrow{\iota_{\lambda,\alpha}^{\mu,\nu}}\mu\otimes \nu\;\;\mbox{with}\;\;
\sum_{\lambda,\alpha} \iota_{\lambda,\alpha}^{\mu,\nu}\circ \pi^{\lambda,\alpha}_{\mu,\nu}=\mbox{id}_{\mu\otimes \nu}\ ,\label{biprodf1} 
\end{equation}
where the summations $(\lambda,\alpha)$ are over $\hat{I}_{\mu,\nu}$ such that $\alpha$ denotes the different copies of 
$\lambda\in \hat{I}_{k,l}$ and for all $(\lambda,\alpha),(\lambda',\alpha')\in \hat{I}_{\mu,\nu}$
\begin{equation}\pi^{\lambda,\alpha}_{\mu,\nu}\circ\iota_{\lambda'\alpha'}^{\mu,\nu}=
\delta_{\lambda,\lambda'}\delta_{\alpha,\alpha'}\,\mbox{id}_{\lambda}\label{biprodf2}\end{equation}
We now introduce the diagrammatic notation for $\iota_{\lambda,\alpha}^{\mu,\nu}$ and $\pi^{\lambda,\alpha}_{\mu,\nu}$ 
by the figures
\begin{equation}\includegraphics[width=5cm]{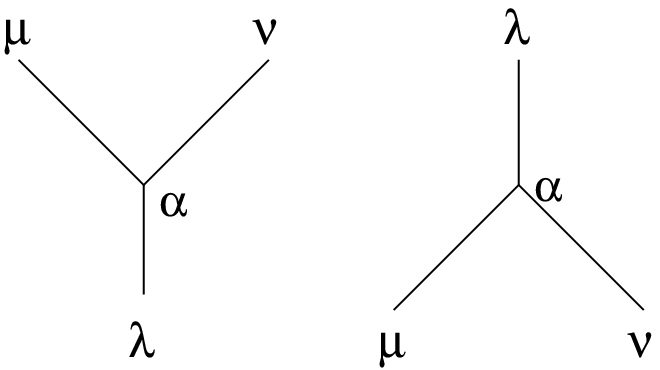}\label{zero}\end{equation}
A line decorated by a simple object $\lambda$ represents the identity id$_{\lambda}$, vertical juxtaposition composition and 
horizontal stacking of diagrams tensor product of arrows.
Then the above compositions (\ref{biprodf1},\ref{biprodf2}) in diagrammatic notation read
\begin{equation}\includegraphics[width=5cm]{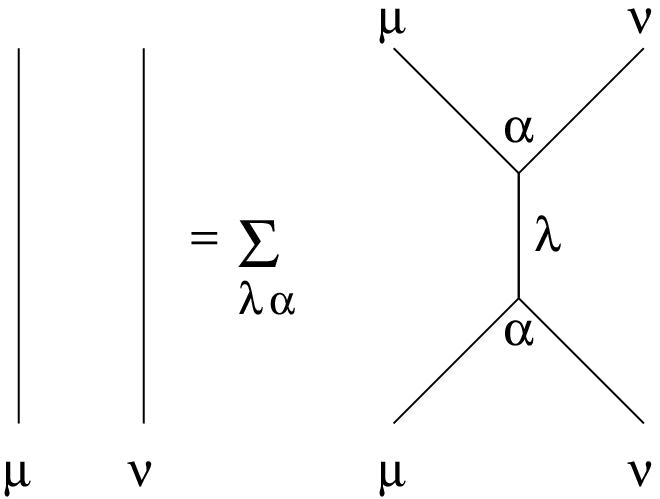}\label{one}\end{equation}
and
\begin{equation}\includegraphics[width=4cm]{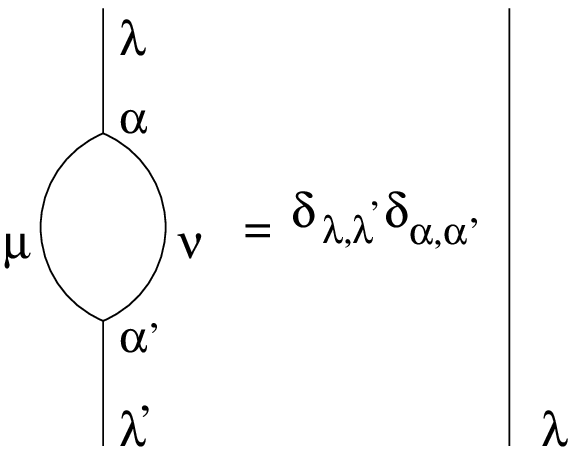}\label{two}\end{equation}
In order to have a diagram calculus in terms of the above, we need the following corollaries of the axioms. Let $a,b$ be objects
in $\cat$ and $m$ be number of simple summands of $a$ and $b$ simultaneously. Then
\begin{equation}\mbox{Hom}_\cat(a,b)\cong\mbox{Hom}_\cat(\oplus_{i\in I} \lambda_i,\oplus_{j\in J} \lambda_j)\cong \times_{i,j}\mbox{Hom}_\cat(\lambda_i,\lambda_j)\cong \times_{i\in I\cap J} 
\mbox{Hom}_\cat (\lambda_i,\lambda_i)\equiv \times_{n=1}^m k\ .\label{il}\end{equation} 
The first isomorphism exists (i.e., the multiset of simples $\{\lambda_i\}$ and $\{\lambda_j\}$ are finite 
due to semisimplicity), the second due to the existence of biproducts 
(the proof is e.g., in Ch. I.17 of \cite{Mitchell}), the third due to the third additional property above 
in the definition of a tensor category.
\lem{For a triple of simples $(\lambda,\mu,\nu)$ in $\cat$ with $N^\lambda_{\mu,\nu}>0$ 
$\{\iota_{\lambda,\alpha}^{\mu,\nu}\}$ and 
$\{\pi^{\lambda,\alpha}_{\mu,\nu}\}$ with $\alpha$ labeling the different copies of $\lambda$ in $\hat{I}_{\mu,\nu}$ 
are basis in Hom$_\cat(\lambda,\mu\otimes \nu)$ and Hom$_\cat(\mu\otimes \nu, \lambda)$, respectively.}
\proof{Linear independence is easily checked using the definition, the spanning property follows from
(\ref{il}), which ensures that the dimension of the Hom spaces coincides with the number of injections/projections. \qed}   
\vspace{0.3cm}
Let us now consider the vector space Hom$_\cat(\lambda,(\mu\otimes \nu)\otimes \rho)$. There is a natural identification
\begin{eqnarray*}
\mbox{Hom}_\cat(\lambda,(\mu\otimes \nu)\otimes\rho)&\to&\oplus_\theta\mbox{Hom}_\cat(\lambda,\theta\otimes \rho)\otimes
\mbox{Hom}_\cat(\theta,\mu\otimes\nu)\\
(\iota_{\theta,\beta}^{\mu,\nu}\otimes\mbox{id}_\rho)\circ\iota_{\lambda,\alpha}^{\theta,\rho}&\mapsto&
\iota_{\theta,\beta}^{\mu,\nu}\otimes\iota_{\lambda,\alpha}^{\theta,\rho}\ ,\end{eqnarray*}
where the direct sum on the rhs. ranges over the multiset 
$\{\theta\in \hat{I}_{\mu,\nu}\,|\,\lambda\in \hat{I}_{\theta,\rho}\}$ and $\alpha, \beta$ denote the basis in the corresponding Homs as
before. We can write a similar line for $\mbox{Hom}_\cat(\lambda,\mu\otimes (\nu\otimes\rho))$.

Although isomorphic, 
$(\mu\otimes \nu)\otimes \rho\neq \mu\otimes (\nu\otimes \rho)$ and there is a natural isomorphism 
\[F^{\mu\nu\rho}_\lambda:\oplus_\theta \mbox{Hom}_{\cal C}(\lambda,\mu\otimes \theta) \otimes(\mbox{Hom}_{\cal C}(\theta,\nu\otimes \rho) 
\to \oplus_\theta\mbox{Hom}_{\cal C}(\lambda,\theta\otimes \rho) \otimes(\mbox{Hom}_{\cal C}(\theta,\mu\otimes \nu)\ .\]
Let us use the basis
$\{\iota_{\theta,\beta}^{\mu,\nu}\otimes\iota_{\lambda,\alpha}^{\theta,\rho}\}_{\theta,\alpha,\beta}$ of the rhs.
and $\{\iota_{\theta,\beta}^{\nu,\rho}\otimes\iota_{\lambda,\alpha}^{\mu,\theta}\}_{\theta,\alpha,\beta}$ of 
the lhs..
We can now write down the above isomorphism in the chosen basis with the help of the diagrams: 
\begin{equation}\includegraphics[width=9cm]{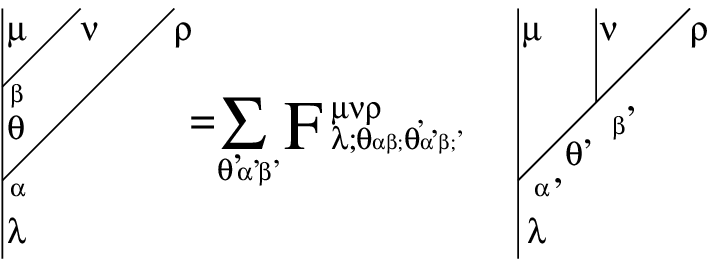}\label{F}\end{equation}
The matrix is called $F$ symbol or (6j) symbol.
If we now consider the vector space Hom$_{\cal C}(d,\mu\otimes \nu\otimes \rho\otimes \theta)$, 
we have 5 different ways to insert brackets 
defining the target object, which leads to the pentagon equation:
\lem{Let us consider the compositions 
Hom$(\lambda, ((\mu\otimes \nu)\otimes \theta)\otimes \rho)\to$ Hom$(\lambda,\mu\otimes (\nu\otimes (\theta\otimes \rho)))$ given
by the figure below with the constituting maps corresponding to the isomorphisms above, $\mu,\nu,\theta,\rho$ are 
labels (of representatives) of simples on top from left to right in each diagram, $\lambda$ is the bottom label. 
Due to the associativity axiom the two compositions are equal, 
which statement schematically reads $FFF=FF$ in terms of the $F$-symbols (the index structure is to be read off from 
the figure). The ``shrinking line'' is thickened in each term of the composition and $r,d$ distinguishes 
the move to the right and down, respectively. 
\begin{center}\includegraphics[width=9.5cm]{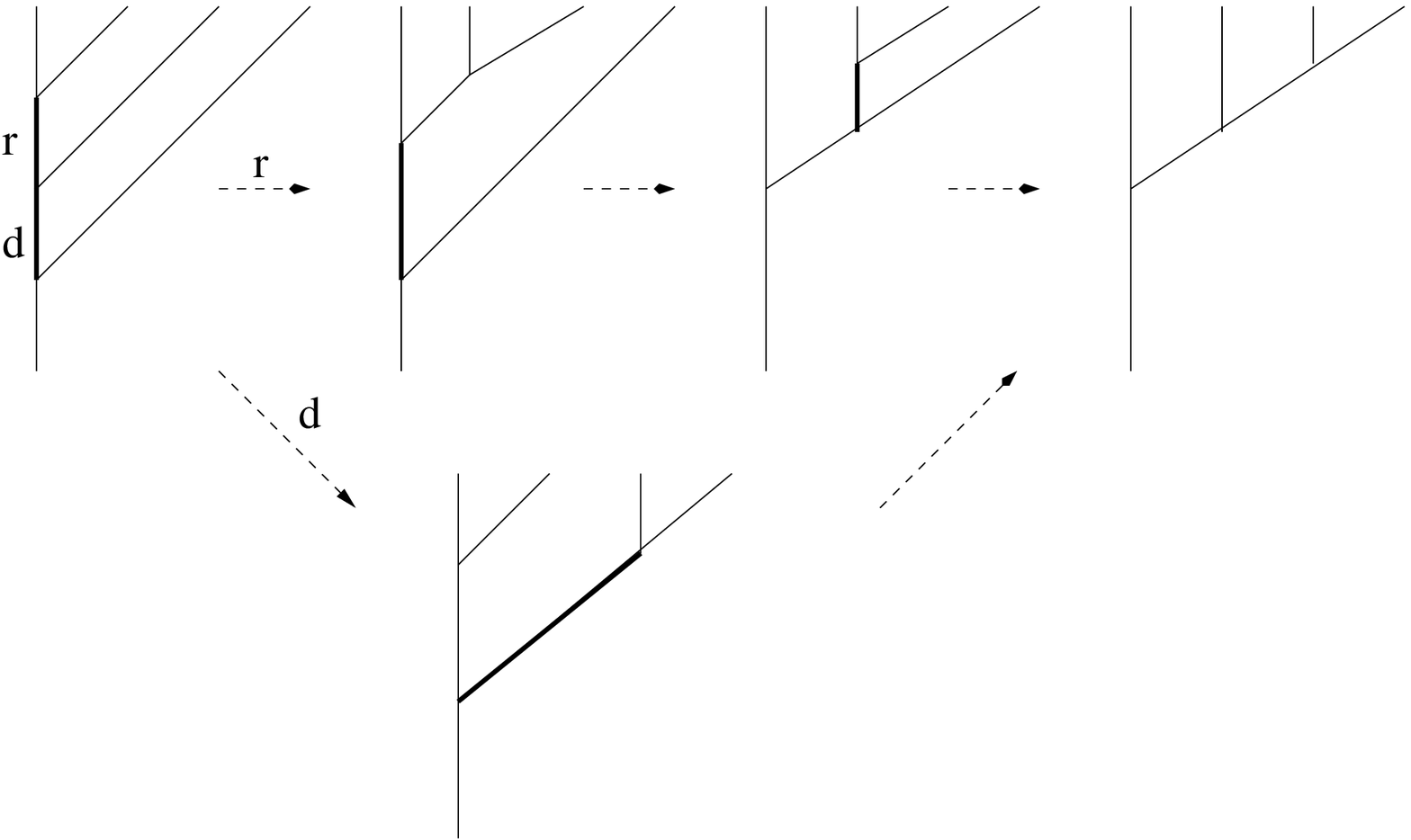}\end{center}
}
Another consequence of the axioms is that for any pair $(\mu,\nu)$: 
$\pi_{\mu,{\bf 1}}^\nu,\,\pi_{{\bf 1},\mu}^\nu,\,\iota^{\mu,{\bf 1}}_\nu,\,\iota^{{\bf 1},\mu}_\nu$ are non-vanishing only if $\mu=\nu$
and one-dimensional (hence the lack of an extra index in the notation). Furthermore, the choices for the injections can be done in 
accordance with (\ref{biprodf1}),(\ref{biprodf2}) such that  
$F^{\mu{\bf 1}\nu}_\rho$ is the basis change from $\{\iota_\nu^{{\bf 1},\nu}\otimes\iota_{\rho,\alpha}^{\mu\nu}\}_\alpha$ to
$\{\iota_\mu^{\mu{\bf 1}}\otimes\iota_{\rho,\alpha}^{\mu,\nu}\}_\alpha$ and a similar statement holds for the projections. 

The totality of the set of labels of simples $I$, finite dimensional $k$-vector spaces $V(\mu,\nu,\rho)\equiv$Hom$(\mu,\nu\otimes \rho)$ 
for all triples $(\mu,\nu,\rho)\in I^3$,the isomorphisms $F^{\mu\nu\rho}_\lambda$ for all tuples $(\mu,\nu,\rho,\lambda)\in I^4$ where
$\lambda$ is a direct summand of $\mu\otimes\nu\otimes\rho$ and a distinguished label ${\bf 1}\in I$ with its properties 
(concerning $V$ and $F$ as in the category) is called a {\em monoidal system}. 
Two monoidal systems $(I,F,V,{\bf 1})$ and $(I',F',V',{\bf 1}')$ are called {\em equivalent} if there is a bijection $I\to I'$ 
mapping ${\bf 1}\to {\bf 1}'$, which induces $V\to V'$ and $F\to F'$ (i.e., e.g., $V(\mu',\nu',\rho')=V'(\mu,\nu,\rho)$).
There is an ``inverse'' procedure to define a tensor category ${\cal C}(I,F,V,{\bf 1})$ from a monoidal system 
up to isomorphism \cite{yama,Kitaev2006}:
{\prop\cite{yama} 
Two tensor categories ${\cal C}(I,F,V,{\bf 1})$ and ${\cal C}(I',F',V',{\bf 1}')$ are isomorphic iff the monoidal systems 
$(I,F,V,{\bf 1})$ and $(I',F',V',{\bf 1}')$ are equivalent.}
\subsection{Diagram calculus\label{dc}}
In the following, we will use diagrams for representing morphisms, which can be composed as horizontal juxtapositions of 
three basic building blocks and vertical concatenation of these whenever the number of lines and their labels match. 
The building blocks are the vertical line labelled $\mu$ representing id$_\mu$ and the trivalent graphs depicted after 
(\ref{biprodf2}). How can we compose two diagrams where target of the first and the source of the second share the same 
labels, but the objects are only isomorphic, not equal? The answer is that the morphisms are vectors in $k$-vector spaces 
and the isomorphisms of objects are also isomorphisms of $k$-vector spaces encoded by the $F$-moves. An example is shown 
in the figure below. We have $f\in$ Hom$_\cat(\kappa\otimes\theta,(\mu\otimes\rho)\otimes \theta)$ and $g\in$ 
Hom$_\cat(\mu\otimes(\rho\otimes \theta),\mu\otimes\nu)$, but they can only be composed as 
$g\circ\alpha_{\mu,\rho,\theta}\circ f$. By composing with the projection decomposition identity (\ref{two}), the diagram is defined
as a linear combination of unambiguous diagrams and the coefficients are the matrix elements of the $F$-moves. 

\begin{center}\includegraphics{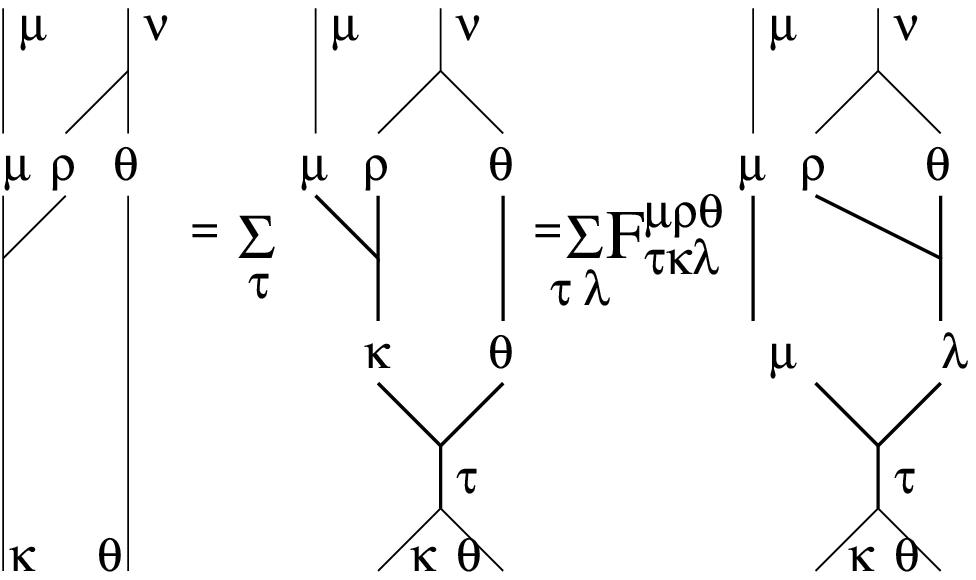}\end{center}
 What if there are more than two vertical lines in a diagram and $\cat$ is not strict? 
To address this question we introduce a ``basis''. A diagram will be defined by the collection of all matrix elements in that basis
regardless of the precise representive of the source and target object in their isomorphism classes. 
Let us choose a tuple of simples $(\lambda_1,\lambda_2,\dots,\lambda_L)$ from $I$ and a finite subset of 
$I_0\subset I$ given. Let $\mu_0\in I_0$.  
We introduce the notation $|\mu\rangle\equiv|\mu_0,\mu_1,\mu_2,\dots,\mu_L\rangle$ 
and $\langle\mu|\equiv\langle\mu_0,\mu_1,\mu_2,\dots,\mu_L|$ for the following
maps in
Hom$_{\cal C}(\mu_L,(\dots ((\mu_0\otimes \lambda_1)\otimes \lambda_2)\dots\lambda_L)$ and 
in Hom$_{\cal C}(\dots((\mu_0\otimes \lambda_1)\otimes \lambda_2)\dots\lambda_L),\mu_L)$, respectively, given by the figure 
below. We shall restrict ourselves to {\em multiplicity-free} tensor category in the following, i.e., 
when dim(Hom$_\cat(\mu,\nu\otimes \rho))\leq 1$ and  dim(Hom$_\cat(\mu\otimes\nu,\rho))\leq 1$ for all triples
$(\mu,\nu,\rho)$ of simples. Otherwise, we would need to label the vertices of the diagrams to get a vector in the corresponding
Hom space.  

\begin{center}\includegraphics[width=8cm]{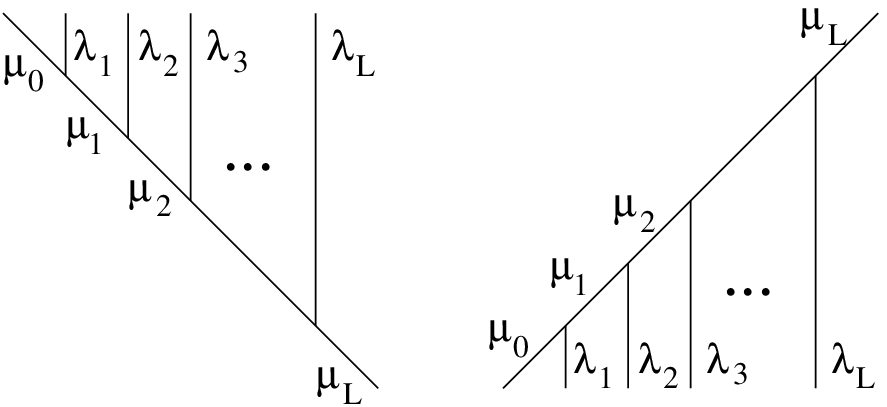}\end{center}
\lem{The set $\{|\mu\rangle\}\equiv\{|\mu_0,\mu_1,\mu_2,\dots,\mu_L\rangle\neq 0\,|\,\mu_0\in I_0,\mu_i\in I$ for $0<i\leq L\}$ 
is finite.}
\proof{The tensor product of simples are isomorphic to the direct sum of finite number of simples, so the range of $\mu_1$ is finite, 
consequently, the range of $\mu_2$ is also finite, and so on. \qed}
\lem{The following orthogonality relations are true 
\[\langle\mu'|\mu\rangle\equiv\langle\mu'|\circ|\mu\rangle=\prod_{i=0}^L\delta_{\mu_i\,\mu'_i}\mbox{id}_{\mu_L}\quad\quad 
\sum_{\{\mu_0\;\mbox{\small fixed}\}}
|\mu\rangle\circ\langle\mu|=\mbox{id}_{\mu_0\otimes\lambda_1\otimes\lambda_2\dots\otimes\lambda_L}\ ,\] where the summation is over
the basis. We shall omit the composition sign
$\circ$ also from expressions like the summands on the rhs.
\proof{For the first equality, we draw the diagram of the composition and recall that $\delta_{\mu_0\,\mu'_0}$ is part of the 
definition of the composition. Then we use the orthogonality relation (\ref{two}) to get rid of
the leftmost bubble and get $\delta_{\mu_1\,\mu_1'}$. 
Then we repeat consecutively the same step proceeding to the right until we get only one line, which represents the 
identity of the object it labels. For the second relation we work towards the opposite direction: $\delta_{\mu_L\,\mu'_L}$ is part
of the definition, then we use the orthogonality relation (\ref{one}) to show that the leftmost middle part of the diagram equals
$id_{\mu_{L-1}\otimes\lambda_L}$, that is two parallel lines with the corresponding labels. Then, we move to the right and 
repeat this step, as a result of which we replace the line labelled by $\mu_{L-1}$ with two lines labelled by $\mu_{L-2}$ and 
$\lambda_{L-1}$
reading from left to right, respectively. 
\qed
\lem{If $I$ is finite and $I_0=I$ then $\phi,\phi'\in $End$_C((\dots(\lambda_1\otimes\lambda_2)\otimes\lambda_3)
\dots\otimes \lambda_L)$ satisfy 
$\phi=\phi'$ iff all their matrix elements in the above basis do.}

\hspace{-0.6cm}We will not prove this here, just mention the two key ingredients: 
(i) semisimplicity ensures the source and target objects are isomorphic
to the direct sum of simple ones and (ii) \cite{fifa} 
for every tensor category ${\cal C}$ there exist a strict tensor category ${\cal C}'=\,$End$\,{\cal C}$,
where the elements of Obj$_{\cal C'}$ are of the form $V':W\to V\otimes W$, see also \cite[E.7.4]{Kitaev2006}. 

We define the equivalence in End$_\cat((\dots(\lambda_1\otimes\lambda_2)\otimes\lambda_3)\dots\otimes\lambda_L)$ as 
$\phi{I_0 \atop =}\phi'$ when all {\em matrix elements} $\langle\mu'|\phi|\mu\rangle$ of $\phi$ defined by
\[\langle\mu'|\phi|\mu\rangle\,\mbox{id}_{\mu_L}\equiv\langle\mu'|\,\circ\,\phi\,
\circ\,|\mu\rangle\ .\]
coincide with those of $\phi'$. Note, that the diagram is well defined only 
for $\mu_0=\mu'_0$, but we extend the definition by defining it to be $0$ for $\mu_0\neq\mu'_0$. 
The figure on the right is the definition of 
$|\mu'\rangle\langle\mu|\in$ Hom$_{\cal_C}((\dots(\mu_0\otimes\lambda_1)\otimes\lambda_2)\otimes\dots\otimes\lambda_L),
((\dots(\mu'_0\otimes\lambda_1)\otimes\lambda_2)\otimes\dots\otimes\lambda_L)$.
\begin{center}\includegraphics[width=10cm]{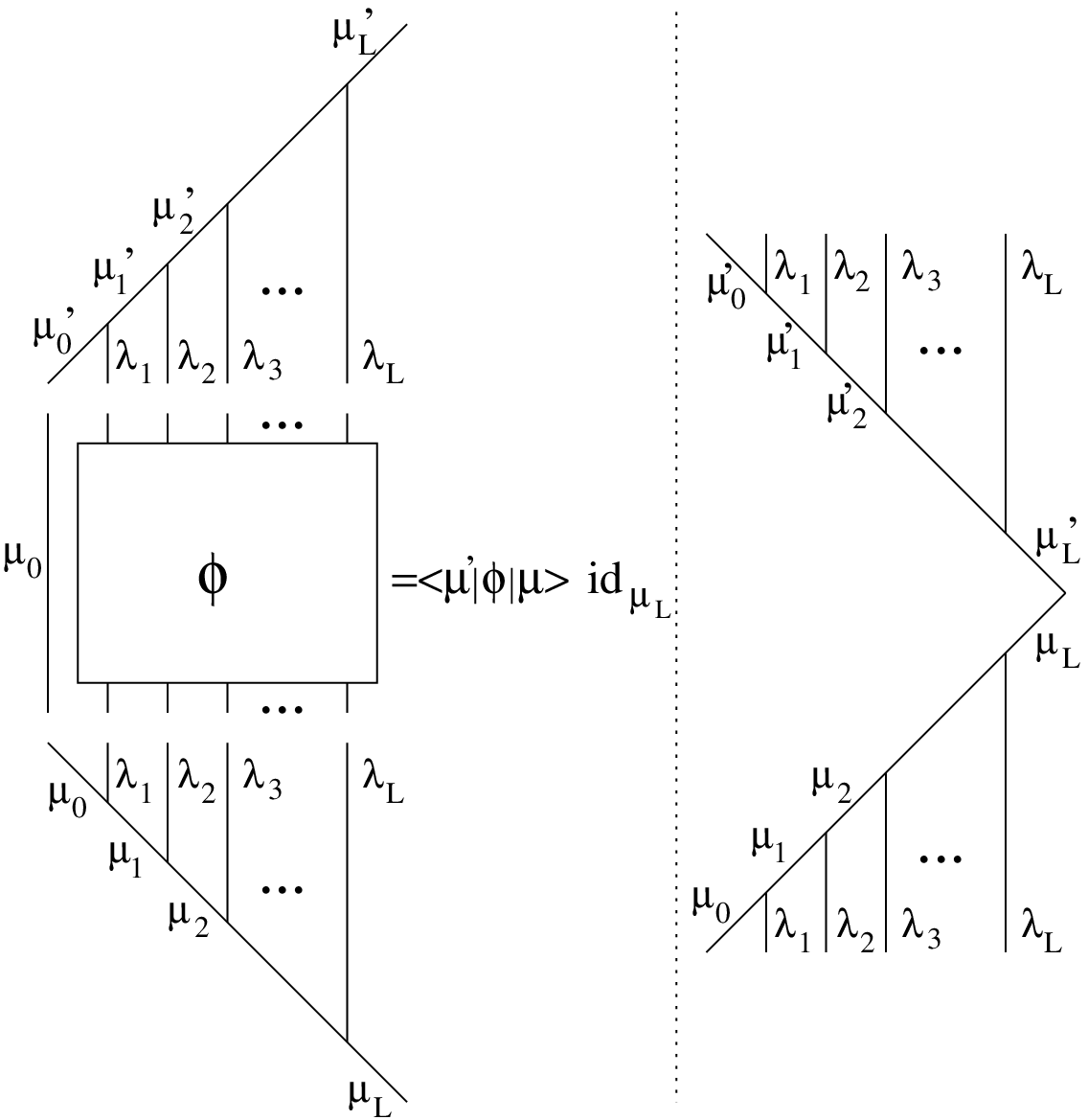}\end{center}
In the following, we will write down expressions like End$_{\cal C}(\mu\otimes\nu\otimes\rho)$ 
even if this defines the source (and target) object up to isomorphism as we always mean the collection of matrix elements in 
the above basis.

\hspace{-0.6cm}{\bf Example:} $U(sl(2))$

The fusion rules of the representation category of the Lie algebra $sl(2)$ is given in the beginning of 
subsection \ref{sc}. Choosing $I_0=\{0\}$, the basis for $L=3$ and $\lambda_1=\lambda_2=\lambda_3=1$ is 
\begin{center}\includegraphics[width=2.7cm]{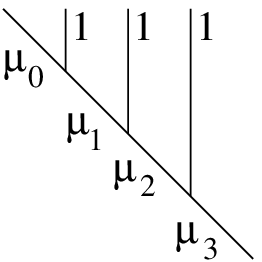}\end{center}
is given by $\{|0,1,0,1\rangle,|0,1,2,1\rangle,|0,1,2,3\rangle\}$. This basis characterizes any 
$\Phi\in$ End$(1\otimes 1\otimes 1)$, that is,
\begin{center}\includegraphics[width=5cm]{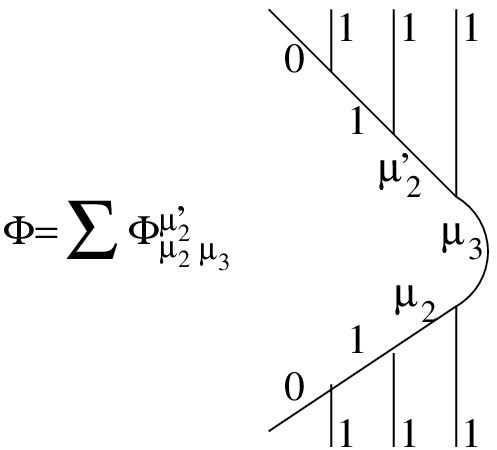}\end{center}
where the only non-zero coefficients are $\Phi_{2,3}^2=\langle 0,1,2,3|\Phi|0,1,2,3\rangle, 
\Phi_{\mu_2,1}^{\mu'_2}=\langle 0,1,\mu_2,1|\Phi|0,1,\mu'_2,1\rangle$ with $\mu_2,\mu'_2\in \{0,2\}$. This
correspond to the fact that $1\otimes 1\otimes 1\cong 0\oplus 0 \oplus 3$, so any element in End$(1\otimes 1\otimes 1)$ 
can be characterized by a 3 x 3 block-diagonal matrix, where the 1-dimensional block is the map $3\to 3$ and the two-dimensional one
is $0\oplus 0\to 0\oplus 0$. 
%
%
%
%
\section{Temperley-Lieb relations} \label{prdr}
Here we state and prove our main theorem. It says that under rather mild conditions, we can define projections 
from the arrow set of a tensor category, which obey relations similar to the Temperley-Lieb algebra. 
Then, demonstrating the equivalence of the monoidal categories and monoidal systems, we 
compute the matrix elements of the projections in the theorem in terms of which the proof in \cite{Finch3} is 
written. 

Let us recall the definition of the Temperley-Lieb algebra $TL_n(\delta)$. Let $k$ be a field, $\delta\in k$ invertible. 
The Temperley-Lieb algebra $TL_n(\delta)$
is generated by the set $k\{U_1,U_2,\dots,U_{L-1}\}$ subject to the relations
\[\begin{array}{ll}U_i^2=\delta U_i&i=1..L-1\\\\
U_iU_j-U_jU_i=0&|i-j|>1\\\\
U_iU_{i+1}U_i=U_i&i=1..L-2\\\\
U_iU_{i-1}U_i=U_i&i=2..L-1\\\\
\end{array}\] 
We will state the main theorem for {\em multiplicity-free} tensor categories, that is, 
dim(Hom$_{\cal C}(\alpha\otimes\beta,\gamma))\leq 1$ for any triple $(\alpha,\beta,\gamma)$ of simple objects. 
{\thm\label{maint} Let ${\cal C}$ be a multiplicity-free tensor category with an algebraically closed base field $k$ 
and \newline 
$\lambda_i,i=1..L,\; \nu_j,j=1..L-1$ simple objects with the property that the spaces 
Hom$_{\cal C}(\lambda_i\otimes \lambda_{i+1},\nu_i)$ are not empty for all $i=1..L-1$. Let us define projections
$p_i=\bigotimes_{i=1}^L\lambda_i\to \bigotimes_{i=1}^L \lambda_i$ by
%
\[\includegraphics[height=3.7cm]{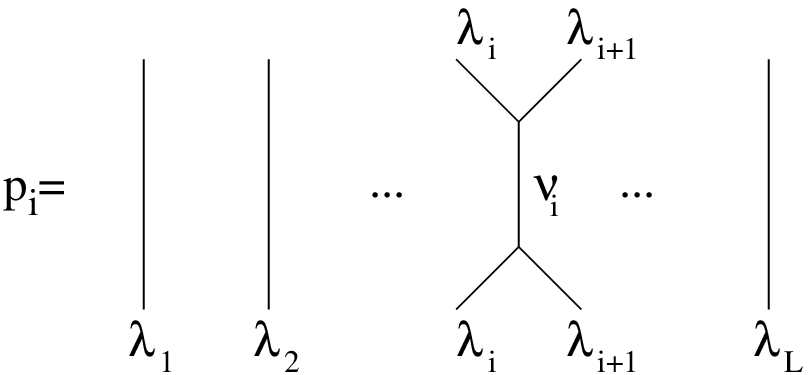}\]
%
Let $\nu_i\otimes \lambda_{i+2}\cong\lambda_i\otimes \nu_{i+1}\cong\mu_i$ for some simple objects $\mu_i, i=1..L-1$. Then 

\hspace{-0.5cm}(i) $p_ip_{i+1}p_i=c_ip_i$ for $i=1,2,\dots,L-2$ and $p_ip_{i-1}p_i=c_{i-1}p_i$ for $i=2,3,\dots L$ with 
\[c_i=\mathlarger{F^{\lambda_i\lambda_{i+1}\lambda_{i+2}}_{\,\mu_i\,\nu_i\;\;\;\nu_{i+1}}\,(F^{-1})^{\lambda_i\lambda_{i+1}\lambda_{i+2}}_{\mu_i\nu_{i+1}\,\nu_i}}\ .\]
(ii) The constant $c_i$ is independent of $i$.}

\hspace{-0.6cm}Note, that in the homogeneous case ($\nu_i=\nu, \lambda_i=\lambda$) the assumption says that $\lambda\otimes\nu\cong\nu\otimes\lambda$ and it is simple. Objects, which satisfy that their fusion with any other simple
object is simple are called simple currents and have been studied in the literature \cite{Schellekens,FSS,Bantay1,Bantay2}. Also note
that for a fixed $i$ the only the simplicity $\nu_i\otimes \lambda_{i+2}$ is needed for the proof of $p_ip_{i+1}p_i=c_ip_i$, whereas
only the simplicity of $\lambda_{i-1}\otimes \nu_i$ is needed for the proof of $p_ip_{i-1}p_i=c_{i-1}p_i$. The fact the
only $\nu_i$ is constrained this way for both relations is crucial for the generalisation for the BMW case \cite{Finch3}.
\proof (i) It suffices to show that the maps given by the diagrams inside the dashed rectangles in the figure satisfy the equality.
\begin{center}\includegraphics[height=4.7cm]{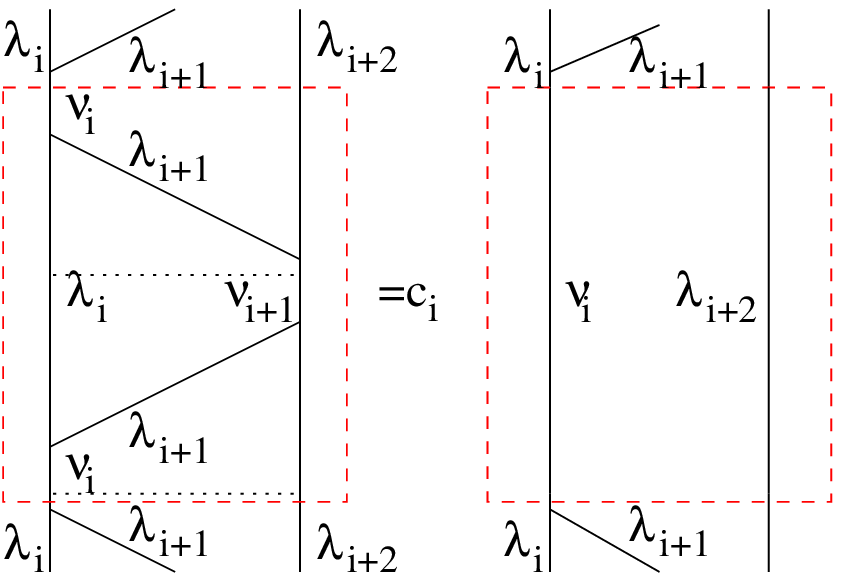}\end{center}
We insert the projection decomposition (\ref{one}) of the identity id$_{\nu_i\otimes\lambda_{i+2}}$ at the lower dotted line and that of 
id$_{\lambda_i\otimes\nu_{i+1}}$ at the upper dotted line and change basis. (The internal lines in of the tree diagrams corresponding 
to base changes are thickened in the figure.)   
\begin{center}\includegraphics[height=6.2cm]{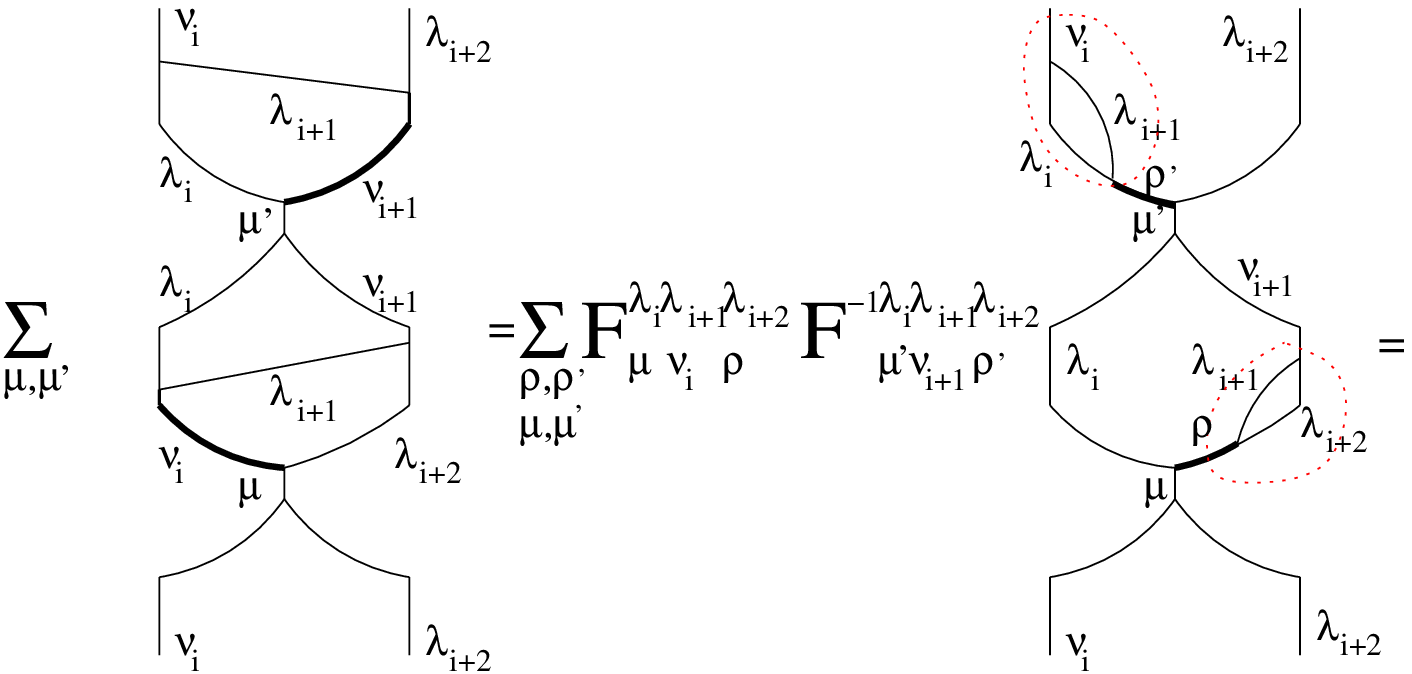}\end{center}
Then, we use relation (\ref{two}) to cancel the summations over $\rho,\rho'$ in the small loops (encircled by red dotted lines 
in the figure) in the rhs.. One also notices that the summations over $\mu$ and $\mu'$ contain only one term $\mu_i$ due to the 
condition of the theorem: 
\begin{center}\includegraphics[height=5.2cm]{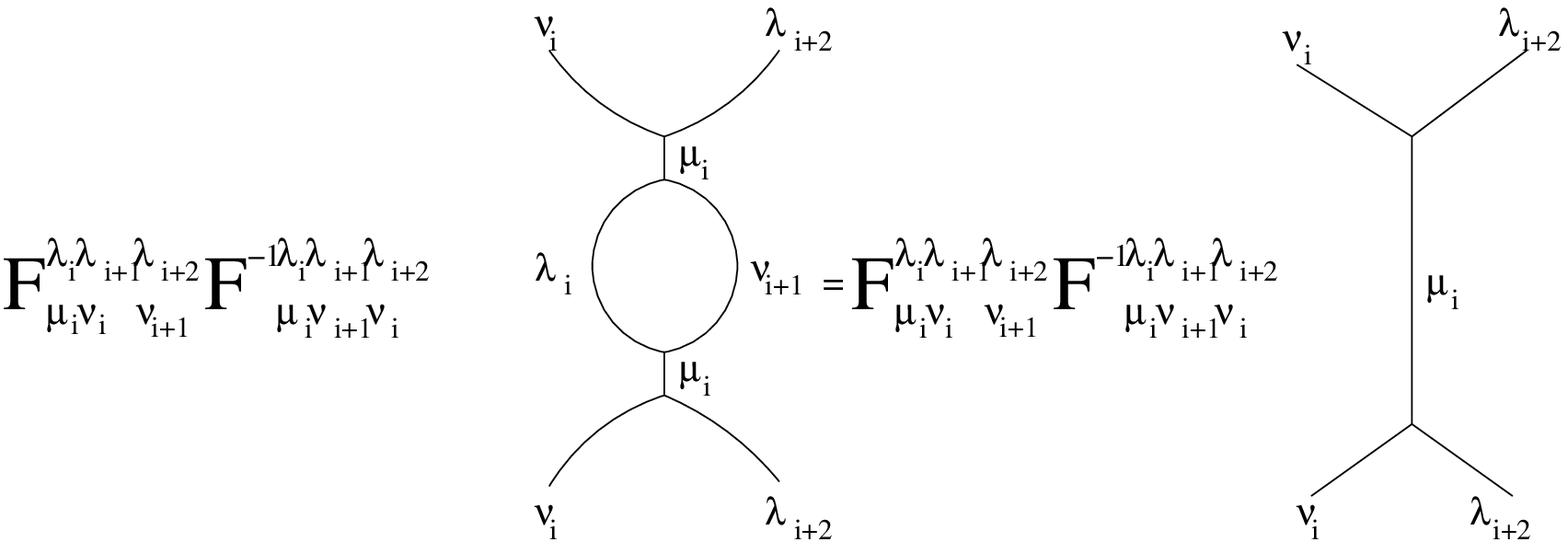}\end{center}
Here the big loop disappeared by (\ref{two}) and the prefactor comes out to be precisely $c_i$. Once again, by the condition 
$\mu_i\otimes\lambda_{i+2}\cong \mu_i$ the diagram is $c_i$ id$_{\nu_i\otimes \lambda_{i+2}}$. 
The proof of the other relation $p_ip_{i-1}p_i=c_{i-1} p_i$ is similar. 

\hspace{-0.6cm}(ii) Due to the assumptions of the theorem $c_i\neq 0$ and it is possible to multiply these projections with scalars ($c_0\in k^*$ arbitrary) as 
\begin{equation}
U_1\equiv c_0 p_1,\quad U_{2i}=\frac{\prod_{j=1}^{i-1}c_{2j}}{c_0\prod_{j=1}^i c_{2j-1}}p_{2i},\quad U_{2i+1}=
\frac{c_0\prod_{j=1}^i c_{2j-1}}{\prod_{j=1}^i c_{2j}}p_{2i+1}\label{diaj}\end{equation}
so that the cubic relations of the Temperley-Lieb algebra 
\begin{equation} U_j U_{j\pm 1}U_j=U_j\label{cubic}\end{equation}
are satisfied for all indices in the defined range. 
However, using these normalisations, and denoting the prefactors by $d_j$, we have 
\begin{equation} U_j^2=d_j U_j \label{loope}\end{equation}
and  we can calculate
\[\begin{array}{c}
(d_{i+2}-d_i)U_i U_{i+2}=d_{i+2}U_{i+2} U_i-d_i U_i U_{i+2}=U^2_{i+2}U_i-U_i^2 U_{i+2}=U_{i+2} U_i U_{i+2}- U_i 
U_{i+2} U_i=\\\\
U_{i+2}U_i U_{i+1} U_i U_{i+2}- U_i U_{i+2} U_{i+1} U_{i+2} U_i=0\ ,\end{array}\]
where we used $U_i U_{i+2}=U_{i+2} U_i$ throughout, (\ref{loope}) in the second, (\ref{cubic}) in (both terms of) 
the fourth equation. Since $U_i U_{i+2}$ is non-zero by construction, the coefficient $d_{i+2}-d_i$ has to vanish. 
Now it is easy to show using the form of $d_i$ from (\ref{diaj}) that $c_i$ is a constant independent of $i$.
\qed}

Some remarks about this result are in order. 

The homogeneity of the constant $c_i$ implies that of $d_{2i}$ and $d_{2i+1}$; the constant $c_0\in k^* $ in (\ref{diaj}) remains
arbitrary. Nevertheless, the statement (ii) of the theorem is a great restriction on the $F$ symbols, which reflects the 
difficulty of finding any inhomogeneous chain with the assumptions of the theorem satisfied.
  
The homogeneous case with rigidity and spherical structure, 
where $\nu={\bf 1}$ the vacuum object is known; since 
$c\equiv F^{\lambda\lambda\lambda}_{\lambda{\bf 1}{\bf1}}
(F^{-1})^{\lambda\lambda\lambda}_{\lambda{\bf 1}{\bf1}} =d^{-2}_\lambda$, 
we have $U_i=p_i/d_\lambda$.

We may consider the periodic case when all indices are defined modulo $L$ and $L$ projections are defined. 
Here for the case of odd $L$ we get the constraint $c_0=\pm 1/\sqrt{c}$, whereas for the case of even $L$ the 
constant remains arbitrary (invertible). 

\section{\label{SW}Schur-Weyl duality}

This subsection is a summary of results from \cite{Paulbook}. 
The construction of the regular representation arising from the work of Andrews Baxter and Forrester 
in the context of the eight-vertex SOS model of statistical mechanics \cite{ABF} will be outlined. The basis, which spans 
the entire multi-matrix structure of the algebra, is given by pairs of sequences, which will be shown to correspond to the diagram 
basis of End$(V^{\otimes L})$, where $V$ denotes the defining two-dimensional irrep. of $U_q(sl_2)$ from now on. 

Define a sequence $\{s\}$ to be an ordered set of positive integers indexed by $i=0,1,2,3,\dots$
\[\{s\}=s_0s_1s_2\dots s_i\dots\]
with the properties $s_0=1$ and $s_i-s_{i-1}=\pm 1$. Note that $s_1=2$ is forced. We will use the notation $\{.\}$ for 
such sequences. If $\{s\}$ and $\{t\}$ have length $n+1$ then $(s,t)$ is an element of $T_n(\delta)$ obtained iteratively as follows:
First for $m$ a positive integer, $m+n$ odd, $m<n+2$ and $p=(n-m+1)/2$ define
\[(e_m,e_m)=((12)^p1m,(12)^p1m)=\left(\prod_{i=1}^p U_{2i-1}\right) E_m^{(2p)}\ ,\]
where $(12)^p$ stands for the subsequence with $12$ repeated $p$ times (e.g., $(12)^3=(121212)$), $1m$ for the subsequence $123\dots m$
and $E_m^{(2t)}$ is an idempotent defined as follows. Let $E_m\in T_n(\delta)$ defined by $E_1=E_2=1$, 
$E_m\in T_{m-2}(\delta)\subset T_n(\delta)$ for $1\leq m\leq n+2$, $E_m^2=E_m$ and for $i=1,2,\dots,m-2$
\[E_mU_i=U_iE_m=0\]
Note that these exist and they are unique. For an integer $t$ $E_m^{(t)}$ is defined in such a way that all indices are translated by $t$ in the 
definition (so $E_m^{(t)}\in T_{m-2+t}(\delta)\subset T_{n+t}(\delta)$).

Returning to the defining defining the POS (pairs of sequences) basis of $T_n(\delta)$ we give an iterative procedure. For $n>1$ let 
$k_n=U_{n-2}(\delta/2)/U_{n-1}(\delta/2)$ with $U_n(x)$ being the $n$-th Chebishev polynomial of the second kind. Suppose 
that the $i$-th element of the sequence $\{s\}$: $s_i=g-1$ a minimum of \{s\}, that is, $s_i=s_{i\pm 1}+1=g$. Then, denoting $s^i$ the
sequence identical to $s$ except $s_i^i=s_i+2$ we have
\[(s^i,t)=\sqrt{k_g k_{g+1}}(1-U_i/k_g)(s,t)\]
and in general $(t,s)=(s,t)^T$ where the latter is obtained by writing the generators of $T_n(\delta)$ in the reversed order.
{\thm \cite{Paulbook} (i) If the operators $(u,s)$ and $(t,v)$ are well defined then:
\[(u,s)(t,v)=\delta_{st}(u,v)\]
(ii) When all defined, the set of operators $(t,v)$ for all pairs $(\{t\},\{v\})$ with sequences of length $n+1$  and equal final entry 
are basis for, and as elementary operators exhibit the entire (multi-matrix) structure of the algebra $T_n(\delta)$.}
\subsection{Correspondence with the diagram bases}
It is a fact that for the quantum integers defined by
\[[n]_q=q^{n-1}+q^{n-3}+\dots+q^{-(n-1)}=\frac{q^n-q^{-n}}{q-q^{-1}}\]
there is a faithful representation  $\theta_q:TL_L([2]_q)\longrightarrow$Aut$(V^{\otimes L})$ of the 
Temperley-Lieb algebra, see e.g., \cite{CFS}. Under this representation the basis $\{$id$,U_i,i=1..L-1\}$ is mapped to
$\{$id$,\underbrace{|\dots |}_{i-1}\,{\cup \atop \cap}\,\underbrace{|\dots |}_{L-i-1},i=1..L\}$ in the given order. The symbols
for the images of $U_i$ mean identity in all tensor factors of $V$ except for the $i,i+1$th: ${\cup \atop \cap}$ is the composition
of the projection to the trivial representation with the injection back to $V\otimes V$.

We can map the basis of $T_L([2]_q)$ in terms of pairs of sequences (POS) from Chapter 6 of \cite{Paulbook} summarized in the
beginning of this chapter to the diagram basis of End$(V^{\otimes L})$ of $U_q(sl_2)$. 
Let $q\in {\mathbb C}, q\neq 0,1$ be fixed. Let us use the notation $\mu-1=\{\mu_0-1,\mu_1-1,\dots,\mu_L-1\}$ for sequences
$\mu$ and denote irreps of $U_q(sl_2)$ by an integer $j: (dim(j)=j+1)$. Define the map $\theta$:
\begin{eqnarray*}
\mbox{TL}_L([2]_q)&\to&\mbox{End}(V^{\otimes L})\\
(\mu,\mu')&\mapsto&|\mu-1\rangle\langle \mu'-1|
\end{eqnarray*}
The two basis are well defined for the same range of the parameters and satisfy ``orthonormality'' 
\begin{eqnarray*}(\alpha,\beta)(\gamma,\kappa)&=&
\delta_{\beta,\gamma} (\alpha,\kappa)\\
|\alpha\rangle\langle\beta\,|\,\gamma\rangle\langle\kappa|&=&\delta_{\beta,\gamma} 
|\alpha\rangle\langle\kappa|\end{eqnarray*}
by construction. The action of the algebra $TL_L([2]_q)$ on the POS basis coincides with its image 
$\theta_q(TL_L([2]_q)$ on the diagram basis:
{\thm Let the matrix elements $U_i$ of the generator of $TL_L([2]_q)$ defined by $U_i(s,t)=\sum_{s'}U_i^{s,s'}(s',t)$. 
Then $\theta_q(U_i)=\sum_s'U_i^{s,s'}|s'-1\rangle\langle t-1|$. That is, the two representations agree.}

Let us look at the basis element 
$=|\mu\rangle\langle\mu'|$ of End$(V^{\otimes L})$. By construction $\mu_L=\mu_L'$ (otherwise $|\mu\rangle$ and}
$\langle\mu'|$ cannot be composed (in the order of writing) and due to the fusion rules of $U_q(sl_2)$ 
 we have that 
$\mu_{l+1}\in \{\mu_l-1,\mu_l+1\}$ if $\mu_l\neq 0$, $\mu_{l+1}=1$ if $\mu_l=0$ 
and similarly for $\mu'_{l+1}$, $l\in \{0,1,\dots, L-1\}$.
\begin{center}\includegraphics{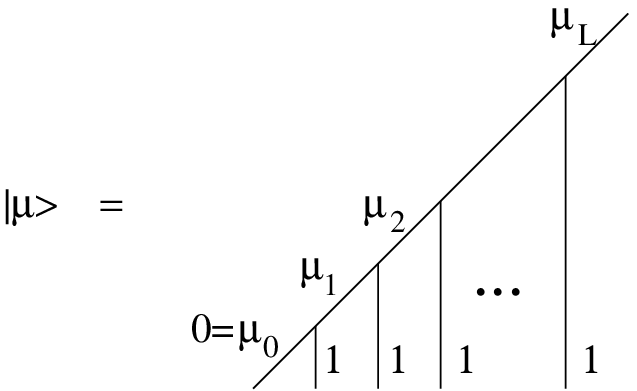}\end{center}
Following the construction in the previous subsection, we look at the diagram corresponding to $(e_m,e_m)$ first. 
For a positive integer $m<L+2$ and $p=(L-m+1)/2$ the image $\theta((e_m,e_m))$ is given by the following composition:
\[\bigotimes_{i=1}^p(V\otimes V)\otimes
(\bigotimes_{j=1}^m V)\xrightarrow{(\otimes_i^p \pi_0)\otimes\,\pi_m} {\bf 0}^{\otimes p}\otimes m
\xrightarrow{\otimes_i^p\iota_0\otimes\,\iota_m}\bigotimes_{i=1}^p(V\otimes V)\otimes(\bigotimes_{j=1}^m V)\ ,\]
where $\pi_n$ is the projection, $\iota_n$ is the injection to the $n$'th irrep, the highest weight module in the $n$-fold tensor
product of $V$. This map in End$(V^{\otimes m})$ satisfies the definition of the unique idempotent $E_m$. 
An example for $L=7, p=2, m=3$ is depicted below. Note, that the location of the trivial irrep. differs from what it is the formula, 
but we are using a fixed basis, in terms of which we can omit or insert vacuum objects without altering the matrix elements.
\begin{center}\includegraphics[width=10cm]{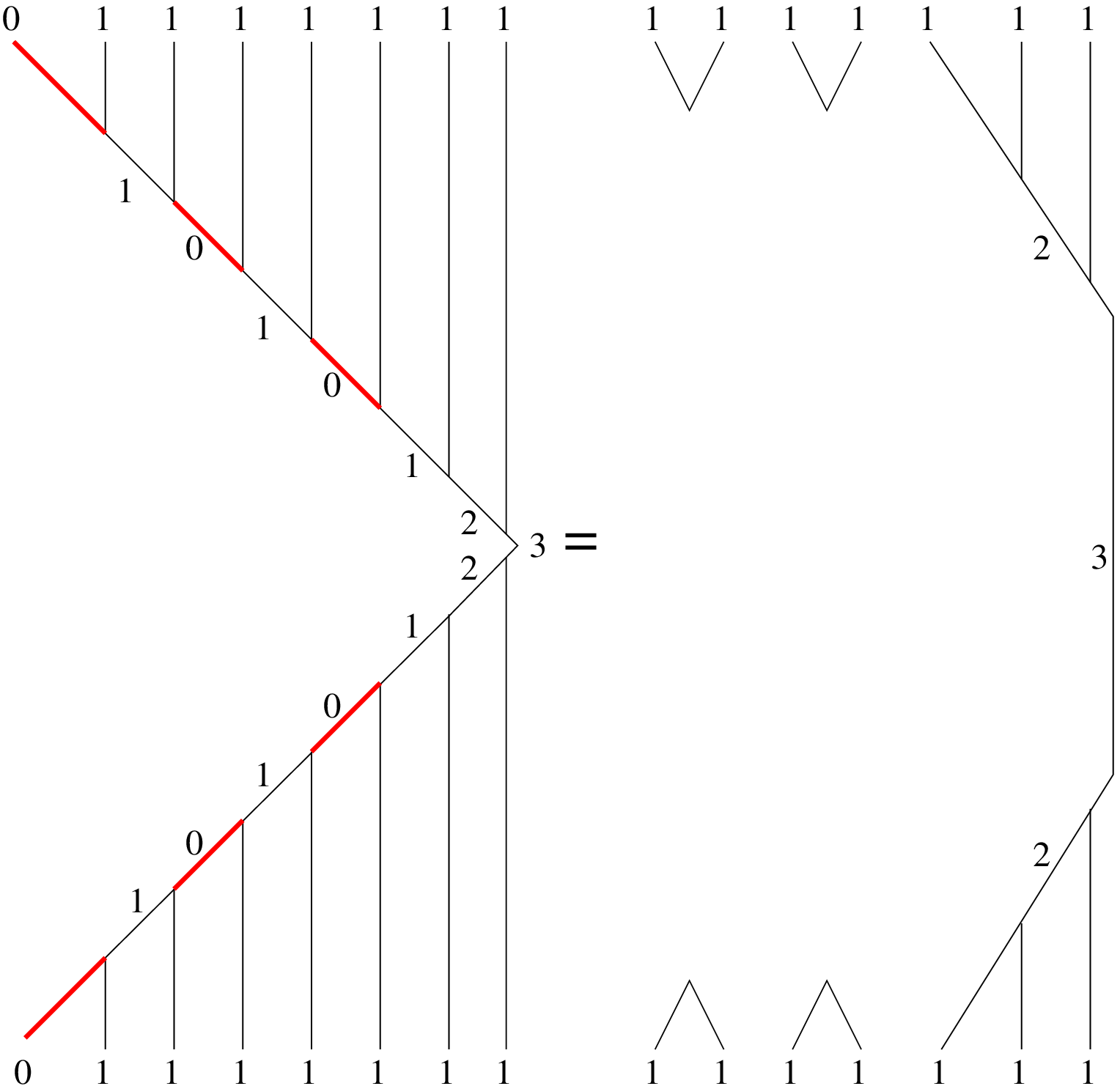}\end{center}
Now let us turn to the proof of the theorem, which follows from the next two lemmas. 
Let us take a basis element with subsequence $\mu_{i-1},\mu_i,\mu_{i+1}$ in the 
diagram notation.
\lem{\label{slope}If $\mu_{i+1}=\mu_i\pm 1=\mu_{i-1}\pm 2$ then $\theta_q(U_i)|\mu\rangle=0$.}  
\proof{Let us denote the loop parameter $\delta\equiv[2]_q$. The proof is the following figure
\begin{center}\includegraphics[width=13cm]{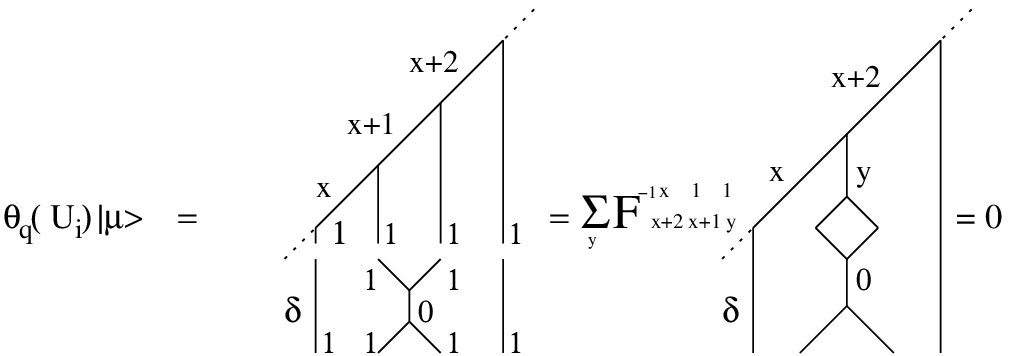}\end{center}
where the last equality holds since, due to (\ref{two}) $y=0$, but the three-valent vertex $(0,x,x+2)$ is not admissible. The statement
for the decreasing sequence is similar. \qed
}
\lem{\label{min}Let $|\mu\rangle$ have a minimum at $i$, i.e., $x\equiv\mu_{i-1}=\mu_i+1=\mu_{i+1}$. 
Let $|\mu^i\rangle$ be the basis element, which is identical to $|\mu\rangle$ except that 
$(\mu^i)_i=\mu_i+2$ (it has a maximum at $i$). Then 
\beq |\mu^i\rangle=\left(\sqrt{k_{x+1}k_{x+2}}-
\sqrt{\frac{k_{x+2}}{k_{x+1}}}\theta_q(U_i)\right)|\mu\rangle\label{v}\eeq
}
The proof can be divided into three steps. 
The state $\theta_q(U_i)|\mu\rangle$ is a linear combination of $|\mu\rangle, |\mu^i\rangle,\newline 
|\dots,x,x-1,x-2,\dots\rangle$, and $|\dots,x,x+1,x+2\dots\rangle$, where $x$ is in the $i$-th position; 
(i) but the coefficients of the latter two states is zero since $U_i^2\sim U_i$ and $\theta(U_i)$ annihilates those 
states. Step (ii) is checking $\langle\mu|\left(\cdot\right)|\mu\rangle=0$, with 
$\left(\cdot\right)$ being the big parenthesis in (\ref{v}).  

Writing explicitly the second term
\bc\includegraphics[width=15cm]{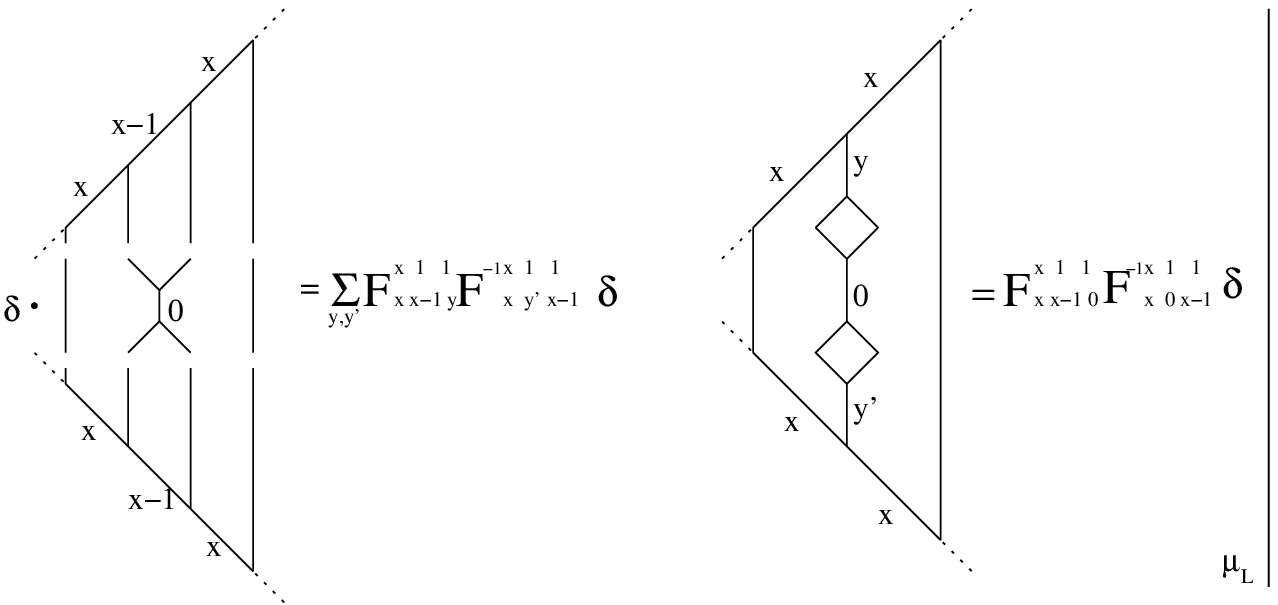}\ec
the statement is equivalent to the equation
\[ F^{x\;\;1\;\;\;\,1}_{x\,x-1\,0}(F^{-1})^{x\,1\;\;\;1}_{x\,0\,x-1}=\frac{k_{x+1}}{\delta}\]
The proof of it for the lowest value of $x=1$ is, for example, a corollary of Theorem \ref{maint}. 
With $\lambda=1$ and $\nu=0$ for the homogeneous case, the lhs. of the above equation is the proportionality factor $c=1/\delta^2$,
which is the correct factor since $k_2=1/\delta$. For the case of general $x$ and step (iii), the normalisation 
one would need to us the explicit construction of 
embedding the irrep. $x$ into $1^{\otimes x}$ as the image of the Jones-Wenzl projection, see e.g., \cite{CFS,ZTQC,KL}.
\section{Other examples}
In case ${\cal C}$ is rigid (i.e., there exist duals and a pair of maps satisfying the rigidity axioms \cite{BK}),
there is a known class of tensor categories, which satisfy the conditions of the main theorem. These are categories, which have among their
simple objects {\em simple currents} $\nu$ for which the map $\nu\otimes\cdot$ is a permutation of the (equivalence classes of) simples. 
There are many rigid tensor categories with simple currents, they correspond to symmetries of the 
fusion algebra and in case the category is also modular they are studied in the context of rational conformal 
theories \cite{Schellekens,FSS,Bantay1,Bantay2}. 

If the category is rigid and pivotal ($V\cong V^{**}$ functorially for $V\in Obj_{\cal C}$) the homogenous case of $\nu$ being a simple 
current reduces to the known case $\nu\cong{\bf 1}$. The explanation is the following. The map 
$Obj_{\cal C}\to Obj_{\cal C}: V\mapsto V\otimes \nu^*$ sends the (isomorpism classes of) simples 
$\nu,\lambda,\mu$ to  ${\bf 1}, \lambda\otimes \nu^*,\mu\otimes \nu^*$, respectively, where, it is easy to show that 
all three are simples. Consequently the representation of the TL algebra arising from the projection in 
End$_{\cal C}(\lambda\otimes \lambda)$ with simple object $\nu$ will be equivalent to that in 
End$_{\cal C}((\lambda\otimes \nu^*)\otimes (\lambda\otimes \nu^*))$ with ${\bf 1}$.

One can study the inhomogeneous case, but the the assertion (ii) of Theorem \ref{maint} limits the possibilities severely. 
For example, if the chain alternates as $(\lambda,\lambda^*,\lambda,\lambda^*\dots)$, we can find some interesting examples e.g., for
$U_q(sl(3))$, where the defining representation is not self-dual.

The most interesting question is to provide a homogeneous example when neither $\lambda$ nor $\nu$ is a simple current or prove that this
cannot occur. This is related to the conjecture in group representation theory, which states that if the tensor product of two irreducible
representations $\rho_1$ and $\rho_2$ of a simple finite group is also irreducible then $\rho_1$ or $\rho_2$ 
is one dimensional. Another issue is multiplicity-freeness, which is a technical asssumption and it may be relaxed. 

Finally we should mention that our result can be extended to a map from braided categories to representations of the Birman Murakami Wenzl 
algebra \cite{Finch3}.

\noindent {\bf Acknowledgements.}
ZK thanks Eric Rowell, Zolt\'an Zimbor\'as, Parsa Bonderson, P\'eter B\'antay and Tobias Hagge for discussions. 
ZK and PM: \epsrc. $\;$ZK also thanks the University of Leeds for support under the Academic Development Fellowship programme.
\appendix
\section{The matrix elements of the projections}
Using only the data of a monoidal system, we can state and prove the theorem. Recall that a  monoidal system is given by the tuple
 $(I,F,V,{\bf 1})$, where
\begin{itemize} 
\item $I$ is a countable set,
\item $V(a,b,c)$ is a finite dimensional k-vector space for all triples 
$(a,b,c)\in I^3$ ($N_{bc}^a\equiv$ dim$(V(a,b,c))$,
\item for all pairs $(a,b)$ $N_{a\bf{1}}^b=N_{{\bf 1}a}^b=\delta_{a,b}$,
\item for each quadrupole $(a,b,c,d)\in I^4$ $\sum_e N_{ec}^d\,N_{ab}^e=\sum_f N_{af}^d N_{bc}^f$ with the range of $e$ ($f$) are defined
by nonzero summand on the left (right) hand side, respectively and an isomorphism
\[ F^{abc}_d:\oplus_f V(a,f,d)\otimes V(b,c,f)\to \oplus_{e\in M} V(e,c,d)\otimes V(a,b,e) \]
\item and for all $a\in I$ basis vectors $l_a\in V(1,a,a)$ and $r_a\in V(a,1,a)$ such that 
$F^{a1b}_c:i_c^{a,b}\otimes l_a\mapsto i_c^{a,b}\otimes r_a$ for any vector $i_c^{a,b}\in V(a,b,c)$.
\end{itemize}
We will write down the matrix elements of the projection. Using only these we can state the theorem without reference to the category, using only
the monoidal system. 
We need to introduce some more structure in order to write down the matrix elements of the projection defined in Chapter \ref{prdr}.
We denote the dual vector space to $V(a,b,c)$ by $V^*(a,b,c)$ and define the product $V(a,b,c)\times V^*(a',b,c)\to \mathbb{R}$:
\[ (e_\alpha,e_\beta^*)\mapsto \langle e_\alpha,e_\beta^*\rangle\equiv\left\{\begin{array}{ll}e_\beta^*(e_\alpha)=\delta_{\alpha,\beta}&a=a'\\
0&a\neq a'\end{array}\right.\]
Note that this product correspond to (\ref{two}) in the category. Now we make can use of the diagram calculus without referring to 
a category. The notion of a diagram is essentially identical to that of Chapter \ref{dc} (we do not call $I$ the set of simple 
objects and we make no reference of $V(a,b,c)$ and $V^*(a,b,c)$ to Hom spaces). 
We fix a basis of $V(a,b,c)$ for all triples $(a,b,c)\in I^3$ to decorate vertexes, the dual basis will be used to
decorate the dual vector spaces. The following two modifications are postulated to lead to equivalent diagrams.   
\begin{itemize}
\item If there is a bubble like in (\ref{two}) we can replace it with a straight line and multiply the diagram with the above 
product of basis elements decorating the two vertices. 
\item If there is a subgraph like in (\ref{one}) with identical label structure (same labels for the upper and lower edges) 
with a summation over all labels of the middle line such that the vector space 
corresponding to the above vertex is not $0$ and another summation over the basis of the vector spaces, then the subgraph can be 
replaced by two straight lines (as on the rhs. of (\ref{one})) with the same labels as the upper and lower edges had.  
\end{itemize}
Finally, using the introduced 
duals we can write down the ``recoupling rule'' corresponding to the upside down of (\ref{F}).
\lem{There is an isomorphism $F'^{abc}_d : \oplus_f V^*(d,a,f)\otimes V^*(f,b,c)\to \oplus_e V^*(d,e,c)\otimes V^*(e,a,b)$ and it 
is given by 
$(F^{-1})^{T\,abc}_{\;\;\;d}$.}
\proof{The equation in diagrammatic form is depicted below. Composing it from below with the upside down of the rhs. 
($\in V(d,a,f')\otimes V(f',b,c))$, we get for the rhs. by evaluating the products in the middle $F'^{abc}_{\,\,def'}$. 
\begin{center}\includegraphics[width=6cm]{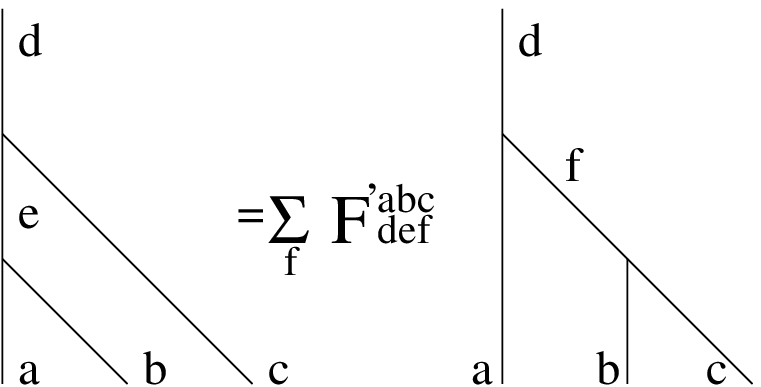}\end{center}
For the lhs. we use the isomorphism $(F^{-1})^{abc}_d$:
\begin{center}\includegraphics[width=10cm]{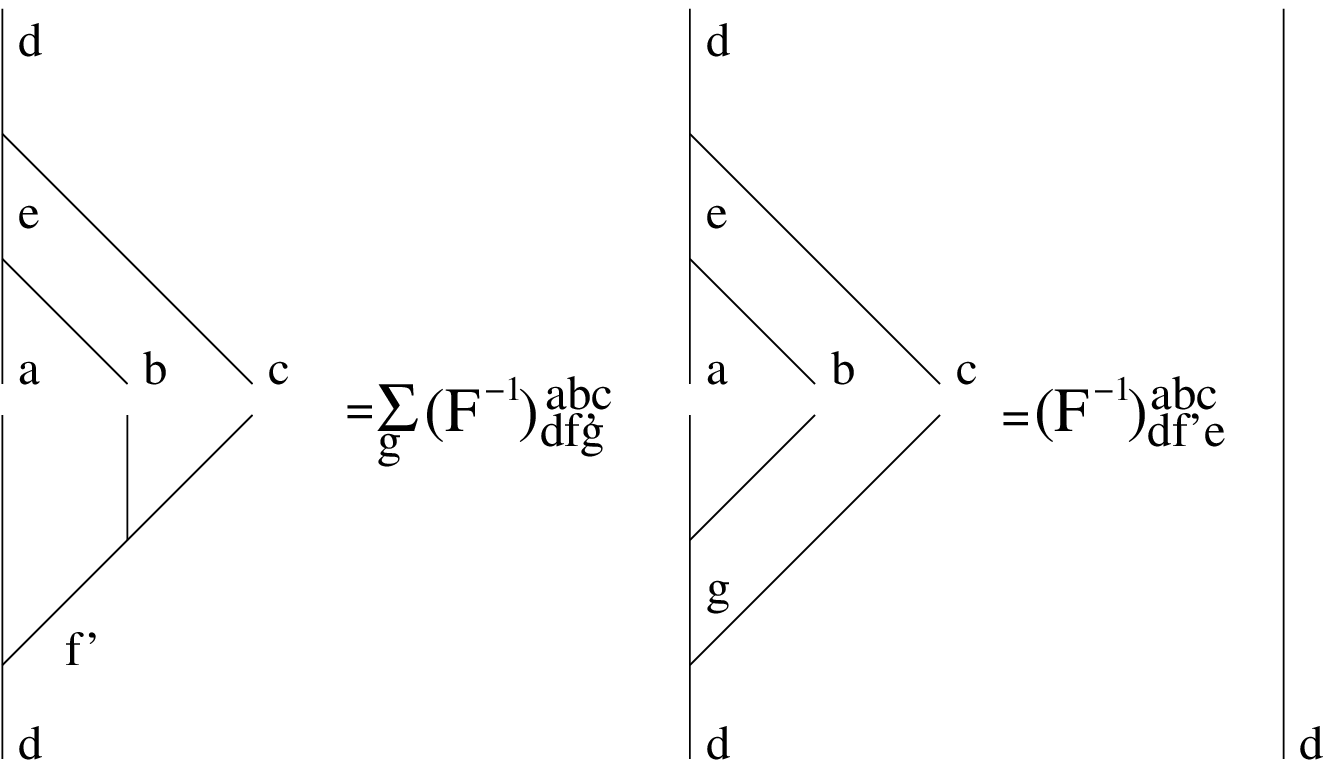}\end{center}
and then evaluate products for the bubbles to arrive at the last equality.
\qed
\vspace{0.3cm}

We can now calculate the matrix of the projection $p_i$ in the defined basis. First we manipulate a part of 
the diagram performing two changes of basis inside as the figure shows. 
\begin{center}\includegraphics[width=14cm]{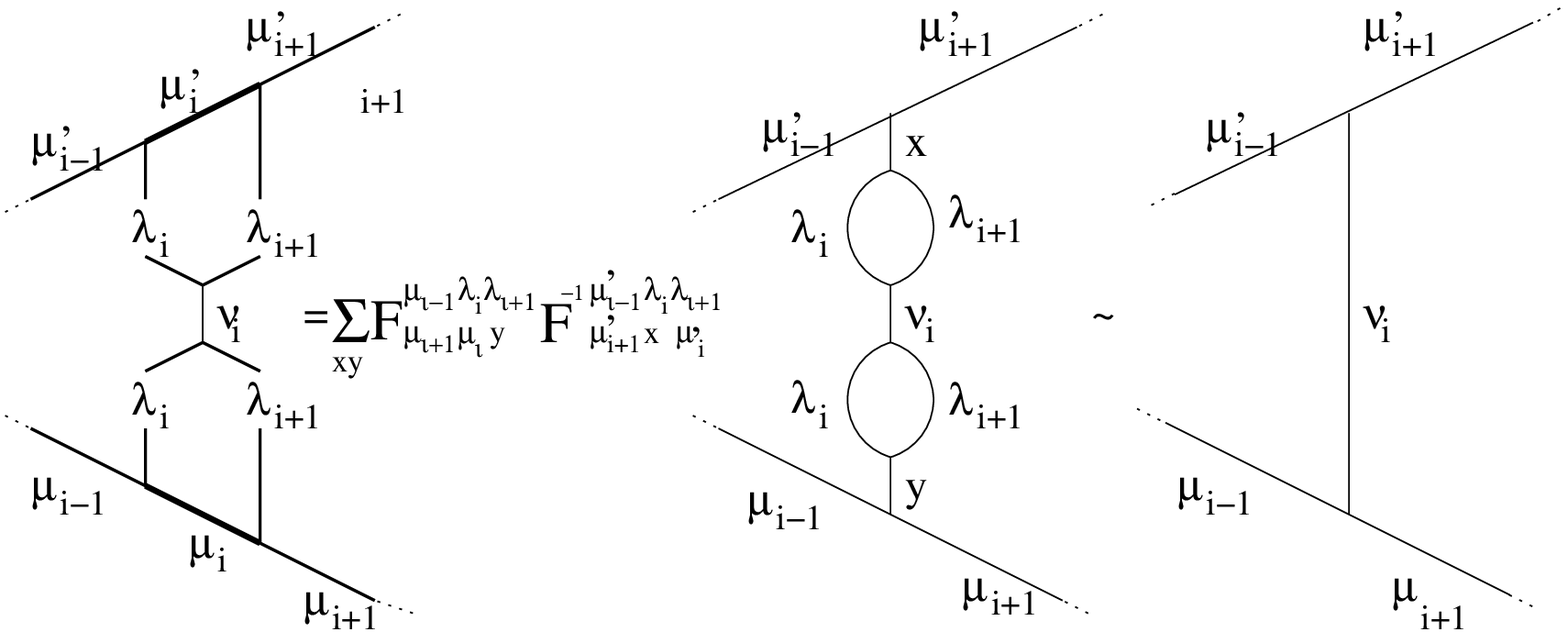}\end{center}
In the (first) equality, we substituted the previously proved result for the isomorphism $F'$. The proportionality is due to the
orthogonality (\ref{two}) applied for both the bubbles in the second diagram. The proportionality factor is 
$F^{\mu_{i-1}\lambda_i\lambda_{i+1}}_{\mu_{i+1}\mu_i\nu_i}(F^{-1})^{\mu'_{i-1}\lambda_i\lambda_{i+1}}_{\mu'_{i+1}\nu_i\mu'_i}$.
Now notice that the above implies 
\[\begin{array}{l}\langle\mu'|p_i|\mu\rangle=F^{\mu_{i-1}\lambda_i\lambda_{i+1}}_{\mu_{i+1}\mu_i\nu_i}
                             (F^{-1})^{\mu'_{i-1}\lambda_i\lambda_{i+1}}_{\mu'_{i+1}\nu_i\mu'_i}\cdot\\\\
\langle\mu'_0\mu'_1\dots\mu'_{i-1}\mu'_{i+1}\dots\mu_L|\mbox{id}_{\lambda_1\otimes\lambda_2\otimes\dots\otimes\lambda_{i-1}
\otimes\nu_i\otimes\lambda_{i+2}\otimes\dots\otimes\lambda_L}|\mu_0\mu_1\dots\mu_{i-1}\mu_{i+1}\dots\mu_L\rangle\\\\
=F^{\mu_{i-1}\lambda_i\lambda_{i+1}}_{\mu_{i+1}\mu_i\nu_i}
                             (F^{-1})^{\mu'_{i-1}\lambda_i\lambda_{i+1}}_{\mu'_{i+1}\nu_i\mu'_i}\prod_{j\neq i}\delta_{\mu_j,\mu'_j}
\end{array}\]
Note that
since $N_{\lambda_i\lambda_{i+1}}^{\nu_i}\neq 0$ was assumed, there are always nonzero matrix elements of $p_i$ provided $I_0=I$.

{}

\end{document}